\documentclass{article}

\usepackage{amsmath,amsfonts,amsthm,amssymb,bm}
\usepackage[margin=1in]{geometry}
\usepackage{hyperref,url}
\usepackage{color}
\usepackage{authblk}
\usepackage{multirow}
\usepackage{graphicx}
\usepackage[table]{xcolor}
\usepackage{array}
\usepackage{booktabs}
\usepackage{caption}
\usepackage{subcaption}
\usepackage{titlesec}
\RequirePackage{nicefrac}
\RequirePackage[textsize=tiny]{todonotes}
\RequirePackage{bm}

\graphicspath{{figs/}}

\newcolumntype{C}[1]{>{\centering\arraybackslash}p{#1}}
\definecolor{EB}{rgb}{0.45,0.67,0.70}
\definecolor{gainsboro}{rgb}{0.86,0.86,0.86}

\usepackage{titlesec}
\titleformat{\paragraph}
{\normalfont\normalsize\bfseries}{\theparagraph}{1em}{}
\titlespacing*{\paragraph}
{0pt}{3.25ex plus 1ex minus .2ex}{1.5ex plus .2ex}

\begin{document}

\title{A rapid and automated computational approach to the design of multistable soft actuators}
\author{Mehran Mirramezani$^*$}
\author{Deniz Oktay}
\author{Ryan P. Adams}

\affil{\textit{Department of Computer Science}, \textit{Princeton University}, \textit{Princeton, NJ, USA 08540}\\
\texttt{\{mehranmir,doktay,rpa\}@princeton.edu}
}
%\texttt{\{mehranmir,doktay,rpa\}@princeton.edu}
%\def\thefootnote{$\dagger$}\footnotetext{Equal contributions}
\def\thefootnote{*}\footnotetext{Corresponding author}
\date{}

\maketitle
\begin{abstract}

We develop an automated computational modeling framework for rapid gradient-based design of multistable soft mechanical structures composed of non-identical bistable unit cells with appropriate geometric parameterization. This framework includes a custom isogeometric analysis-based continuum mechanics solver that is robust and end-to-end differentiable, which enables geometric and material optimization to achieve a desired multistability pattern. We apply this numerical modeling approach in two dimensions to design a variety of multistable structures, accounting for various geometric and material constraints. Our framework demonstrates consistent agreement with experimental results, and robust performance in designing for multistability, which facilities soft actuator design with high precision and reliability.
%This can make the design loop tighter and make it easier to create such multistable mechanical structures with high-speed high-precision soft actuation capability.

\vspace{3mm}

\end{abstract}
\section{Introduction} \label{sec:Intro}
Soft robots, which are comprised of compliant structures instead of rigid bodies and links, enable a wide range of smooth and complicated motions and functionalities.
In particular, soft robots facilitate safe and adaptive interactions with humans and complex environments~\cite{BeckerEtal2017,WangEtal2020}, biomedical device design~\cite{PayneEtal2017}, microsurgical instrumentation~\cite{ZanatyEtal2019}, manipulation of delicate objects~\cite{DeimelBrock2015}, navigation through confined space~\cite{HawkesEtal2017}, etc.
Despite several advantages, the structural compliance of soft robots makes it challenging to achieve high-speed performance, high-precision motion control, and high-strength force output.
Mechanical instability such as snap-through buckling is one approach to tackling these challenges through designing ``tunable'' multistable structures, enabling rapid movements and amplified force output when snapping from one stable configuration to another.

Multistability is an increasingly important aspect of compliant mechanical structures, allowing reconfigurable systems that are locally stable and so can stay in position without energy consumption~\cite{ShanEtal2015,ZhangEtal2016}. 
%This offers several unique soft robotic applications in acoustic~\cite{WangEtal2014,FrenzelEtal2016} and photonic~\cite{LiEtal2012} bandgap control, high energy absorption and release~\cite{CorreaEtal2015,SunEtal2019}, and mechanical actuation~\cite{JeongEtal2019}. 
%For example, in motion control, multistability enables large-magnitude movements in a pre-designed way by switching between different stable configurations, making them work well even with open loop control~\cite{WingertEtal2002,HuEtal2015}.
%This can reduce the number of parts, weight, and assembly time of soft actuators while increasing precision and reliability in soft robotic applications~\cite{SchaffnerEtal2018,GorissenEtal2020}.
Multistable structures are generally designed by combining several bistable structures with snap-through buckling~\cite{SimitsesEtal2006}, a kind of elastic instability in which a structure jumps from one stable configuration to another when external stimuli beyond critical values are applied.
Initially curved beams are commonly adapted as basic elements to create such snap-through induced multistable structures (SIMS)~\cite{OhKota2009,CheEtal2016,HuEtal2021}, (cf.\ Figure~\ref{fig:ms_f1}).
In addition to ease of manufacturing via additive manufacturing, these multistable structures demonstrated promising soft robotic applications in high energy absorption and release~\cite{CorreaEtal2015,SunEtal2019}, fast object grasping ~\cite{BaumgartnerEtal2020}, and rapid/precise locomotion control in a pre-designed way by switching between stable configurations~\cite{WingertEtal2002,HuEtal2015}.
Fundamentally, these locally stable configurations (and the transitions between them) are determined via complex geometries~\cite{IonEtal2017,IonBaudisch2019} that, when combined with appropriate material choices, take advantage of highly nonlinear snap-through buckling~\cite{ZhangEtal2020}.

Unfortunately, this nonlinearity and geometric complexity means that current approaches to multistablity design involve a combination of deep human expertise and experimental trial-and-error~\cite{JeongEtal2019,VangelatosEtal2020,IonEtal2017,IonBaudisch2019}, leading to slow iteration times and suboptimal designs.
Simulations via, e.g., the finite element method (FEM)~\cite{Hughes2012}, can in principle enable comprehensive analysis of such structures, but complex geometries and large deformations require infeasibly costly high-resolution simulations that are prone to numerical instabilities and can be sensitive to method parameters. 
These computational challenges and more importantly non-differentiable nature of most FEM solvers limited their use as a tool for optimization and inverse design of such structures, which requires computation of gradients of an objective with respect to design parameters.
Althought adjoint methods~\cite{CaoEtal2003} are used as computationally effective approaches to numerically compute gradients for such problems~\cite{MinEtal2015,BrucknerEtal2017,ColburnEtal2021}, they require manually deriving derivatives for each input parameters and handling nonlinearities, discontinuities, constraints and boundary conditions, making their implementation non-trivial in many problems.
Therefore, FEM analysis are often employed to validate experimental observations or to understand multistable structures responses under different loading conditions or changes in geometric parameters~\cite{HuEtal2021,ZhangEtal2021}.

Theoretical one-dimensional elastic beam models are widely utilized to mitigate theses computational cost and challenges to inform designs in a more timely manner~\cite{BrennerEtal2003, CheEtal2016,MeaudChe2017,YangMa2020}. For example, Brenner et al.,~\cite{BrennerEtal2003} used a linear beam model and an adjoint optimization technique to design a bistable switch with optimal force ratio.
These models can reasonably capture qualitative behavior of multistable structures e.g., sequence of deformations and multistability~\cite{CheEtal2016}, however, their quantitative predictions such as force amplification magnitude and stability locations could be inaccurate in many cases for a couple of reasons:
First, the simplified assumptions of beam theories such as small deformation are not valid in multistability analysis with large deformations.
Second, one-dimensional beam theories are based on linear-elastic constitutive models, however, soft actuators are mainly fabricated by hyperelastic materials to leverage their extensive elastic deformation and energy absorption capability. %\todo{RPA: If we're going to say that other people's work is inaccuate, we need to bring the receipts.}
Therefore, a rapid differentiable FEM-based solver for high-resolution nonlinear large-deformation modeling of hyperelastic material can significantly improve accurate optimization and inverse design strategies of soft multistable actuators.  
  
Here we propose an automated framework that enables the rapid design of SIMS, given only a desired pattern of multistability described by a strain energy-displacement curve.
The overall approach to this inverse design problem is to rapidly optimize the geometry and material property of a mechanical structure to induce a target strain energy-displacement curve, which is widely used to perform multistability analysis where its local minima determine stable configurations of the structure.
First, we construct a parametric space of mechanical structures via programmatically-defined unit cells.
Each unit cell of the structure has parameters determining local deformation behavior e.g., geometric parameters of a curved beam, enabling highly nonlinear macroscopic properties. 
We develop a differentiable solver based on isogeometric analysis (IGA)~\cite{HughesEtal2005} that enables simulation of the statics of large-deformation continuum mechanics, while also providing the Jacobian of computed outputs with respect to input parameters and internal variables, making it possible to perform efficient optimization.
We take advantage of IGA nature representing both geometry and solution in the same basis functions, which allows us to directly create a differentiable map from geometry parameters to solution basis, and also to a specified loss function (defined by the solution) using adjoint methods.
Adjoint methods~\cite{CaoEtal2003} along with automatic differentiation~\cite{BaydinEtal2018} are utilized in our framework to efficiently compute the gradients, which also alleviate customized algorithmic and implementational overhead challenges of adjoint methods.
While automatic differentiation has been playing critical role in machine learning, its potential applications in various engineering problem-solving techniques have drawn more interest in recent years~\cite{LindsayEtal2021,VigliottiEtal2021}.

We leveraged the native differentiability and GPU acceleration of JAX~\cite{JAX}, the successor to the Autograd tool that has been widely used in machine learning.
Recently, several differentiable physics-based simulators in JAX have shown promise in other scientific computing fields, e.g., molecular dynamic simulation~\cite{SchoenholzCubuk2020}, thermodynamics modeling~\cite{Guan2022}, and computational fluid dynamics simulations~\cite{KochkovEtal2021,BezginEtal2023}.
The end-to-end differentiability of our solver enables the use of modern large-scale constrained optimization techniques to directly optimize geometrical and material properties of a model to achieve desired functionality.
We additionally fabricate prototype designs using additive manufacturing and compare simulation results with experiment.
To the best of our knowledge, this is the only fully automated framework utilizing numerical solution of continuum mechanics for designing snap-through multistable structures with desired multistability behavior, which can make the design loop tighter and make it easier to create such structures with high-speed high-precision soft actuation capability.
%\todo[inline]{RPA: Needs a related work section.}

\section{Methods}\label{sec:Methods}
\subsection{Structural geometry}\label{sub:geom}
\begin{figure}[h!] 
\centering
\includegraphics[width=1.0\textwidth]{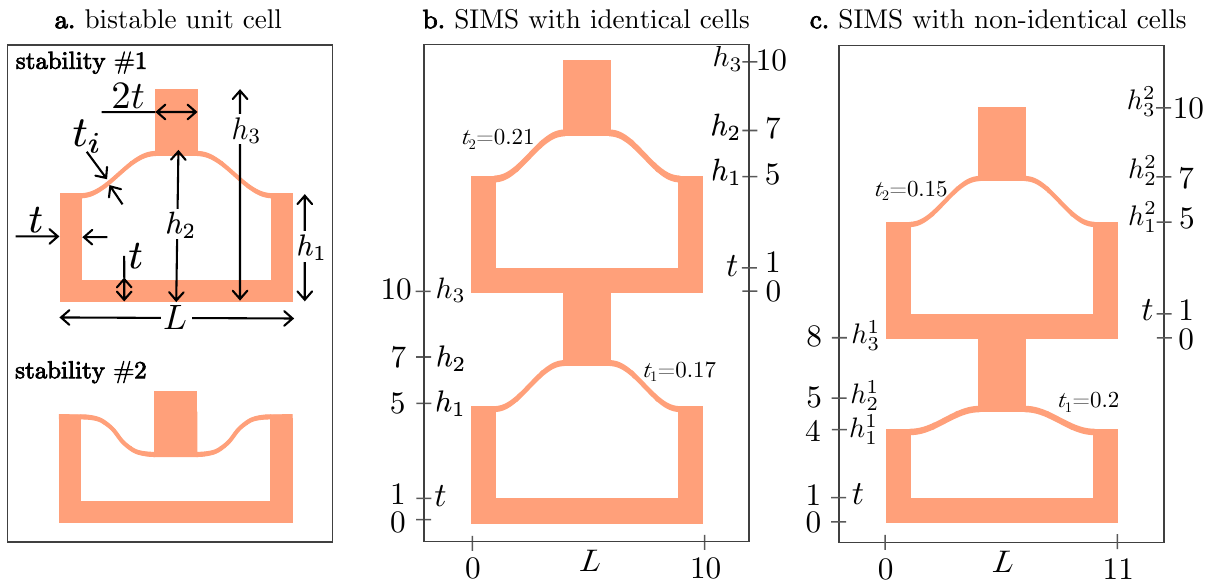}
\caption{\textbf{a)}\ A schematic of a snap-through bistable unit cell with geometric parameters and stable configurations.\ \textbf{b)}\ A representative SIMS constructed from two identical unit cells, but different in-plane beam thickness at each row, with seven geometric design parameters to optimize.\ \textbf{c)}\ A representative SIMS constructed from two non-identical unit cells, with 10 geometric design parameters to optimize.} 
\label{fig:ms_f1}   
\end{figure}

In this study, the design of multistable structure consists of multiple bistable unit cells that each includes a sinusoidal curved beam with snap-through behavior and a frame.
The initial shape of the centerline of a fixed-thickness curved beam $w_0$ can be expressed as:
\begin{equation}
w_0(x) = \frac{h_2-h_1}{2} \bigg [1-\cos{\frac{2 \pi x}{L-4t}} \bigg] \label{beam},
\end{equation}
where $L$ is the length of each unit cell and $t$ is the thickness of the frame. The out-of-plane thickness of the structure is always one cm in our modeling. Other geometric parameters and their descriptions are indicated in Figure~\ref{fig:ms_f1}a for a representative unit cell. The beam's in-plane thickness for each cell is denoted as $t_i$, which is different in each row and considered a tunable geometric parameter to control the sequence of deformation in SIMS with identical unit cells~\cite{CheEtal2016} (Figure~\ref{fig:ms_f1}b). This results in $5+n$ geometric design parameters (GDP) to optimize in our framework to accomplish a desired multistability pattern for a structure with $n$ rows, i.e.\ GDP $= [L,t,h_1,h_2,h_3,t_1,\cdots,t_n]$. Here we also constructed structures with non-identical cells to achieve a broader range of multistability designs, by adding more flexibility to stability locations, that increases the number of design parameters to $2+4n$ where $L$ and $t$ are constant among all the cells (Figure~\ref{fig:ms_f1}c), i.e.\ GDP $= [L,t,h_1^i,h_2^i,h_3^i,t_i]$ for $i = 1,\cdots,n$. 

\subsection{Fabrication and experiments}\label{sub:exp}
To validate our proposed numerical method and study multistability behavior of structures, a series of samples were fabricated via a calibrated Prusa i3 MK3S 3D printer using stl files generated automatically in our computational framework for each model. The printing material is thermoplastic polyurethanes (TPU), known as a material with excellent abrasion resistance and large elastic strain recovery performances~\cite{JaglinskiEtal2006}. Moreover, TPUs have a wide range of Young's modulus ($E$) from 10 to 10000 MPa, which is considered as a parameter in our framework to optimize the design. The value of Young's modulus for each simulation is reported accordingly in different sections. Here, we assumed a fixed Poisson's ratio value ($\nu = 0.46$) based on reported values from prior studies on 3D printed TPU materials~\cite{YangMa2019}. Material properties of the experiment presented in Section~\ref{sub:bistablity} were characterized by standard tensile measurements according to ASTM D638-14~\cite{HaEtal2018}. The measured Young's modulus is approximately $75$ MPa.

% https://www.amazon.com/dp/B07H589LYD/ref=pe_386300_440135490_TE_item
Quasi-static uniaxial loading tests were conducted using a customized compression fixture (shown in Figure~\ref{fig:ms_f2}e) with a 50N tension/compression load cell and a loading rate of 23 mm/min in a displacement-control manner.
Extra parts for connecting samples to the testing system and capturing appropriate fixed bottom and guided side boundary conditions were designed and printed using stiffer carbon-fiber-infused nylon.  

\subsection{Numerical methods}\label{sub:numerical}
\subsubsection{Overview}\label{sec:overview}
We develop a custom isogeometric analysis solver~\cite{HughesEtal2005} for large-deformation plane stress statics problems, with the goal of having a robust and end-to-end differentiable simulator for multistable structure design. Given geometric design parameters, our solver will simulate the multistable structure described by the parameters, and compute the gradient (adjoint) with respect to each design parameter according to a specified downstream loss function; in our case the loss function is $\|\textrm{SE}(D)-\textrm{SE}_t(D)\|_2$, where $\textrm{SE}_t(D)$ and $\textrm{SE}(D)$ are the prescribed (target) and actual strain energy vectors of the structure computed at $N$ displacement loading increments $D=[0,d_1,\cdots,d_{N-1}]$ that are applied at the top edge of a structure. The loss function is then minimized via gradient-based optimization technique to achieve the target strain energy-displacement response.
%\todo{RPA: Is this the norm of a function?  Is it done via integration? Over what range of strains?}

%\begin{figure}
%     \centering
%     \begin{subfigure}[b]{0.48\textwidth}
%         \centering
%         \includegraphics[width=\textwidth]{patch_parent.pdf}
%         \caption{}
%         \label{fig:patchparent}
%     \end{subfigure}
%     \hfill
%     \begin{subfigure}[b]{0.48\textwidth}
%         \centering
%         \includegraphics[width=\textwidth]{constraints.pdf}
%         \caption{}
%         \label{fig:constraints}
%     \end{subfigure}
%        \caption{(a) Visualization of the decomposition of our unit cell into patches. Each patch is pulled back into a parent space, which is a B-spline knot span. Quadrature and integration happens in this space. (b) Visualization of the constraints between control points in different patches, for a representative beam composed of four patches with a fixed top edge. Each control point is labeled with the index it takes in the global degree of freedom array, where each control point corresponds to two degrees of freedom in 2D space. Note that even though there are $64$ control points in the local array (with $128$ degrees of freedom), there are only $27$ control points in the global array (with $54$ degrees of freedom) due to incidence constraints, Dirichlet boundary conditions, and rigidity constraints.}
%        \label{fig:numerical-methods}
%\end{figure}

\begin{figure}[h!] 
\centering
\includegraphics[width=1.0\textwidth]{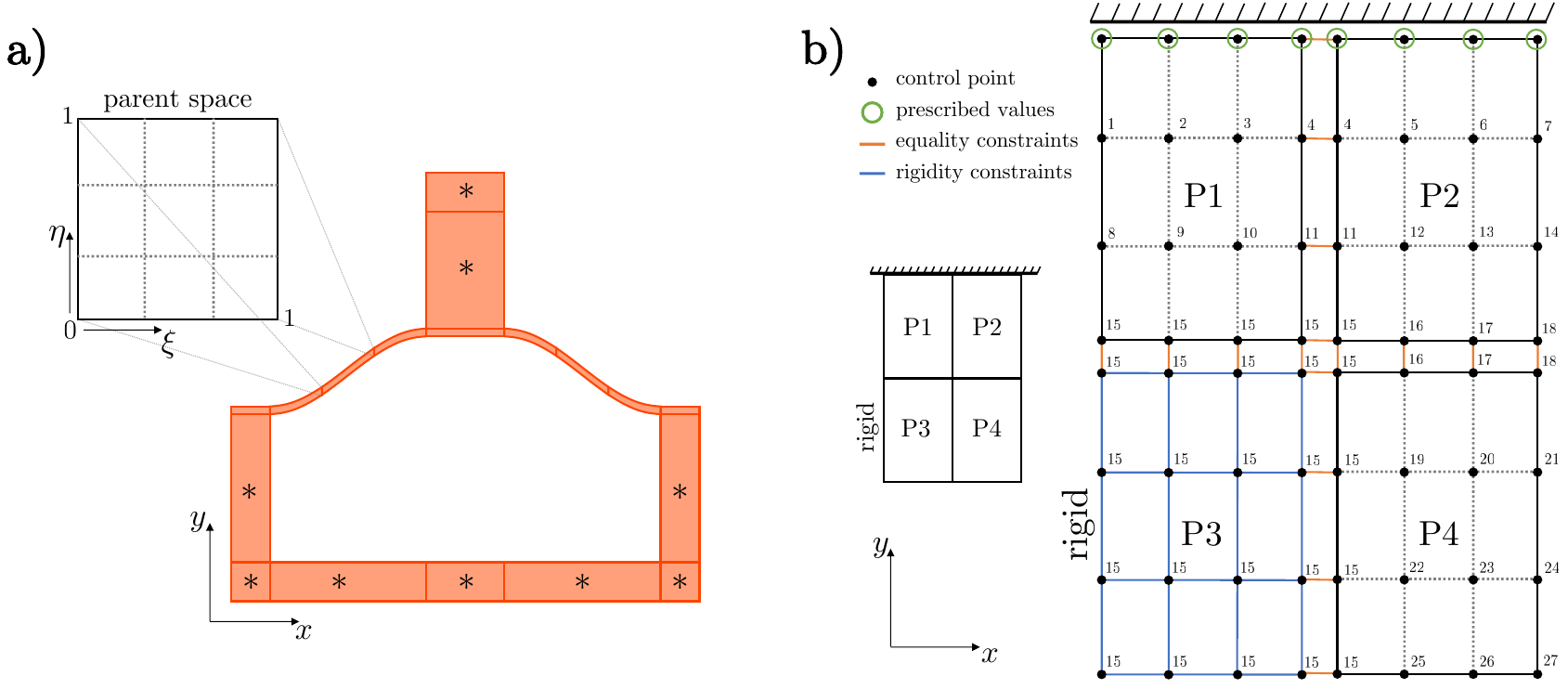}
\caption{\textbf{a)}\ Visualization of the decomposition of our unit cell into patches (rigid patches are starred). Each patch is pulled back into a parent space, which is a B-spline knot span. Quadrature and integration happens in this space.\ \textbf{b)}\ Visualization of the constraints between control points in different patches, for a representative beam composed of four patches with a fixed top edge. Each control point is labeled with the index it takes in the global degree of freedom array, where each control point corresponds to two degrees of freedom in 2D space. Note that even though there are $64$ control points in the local array (with $128$ degrees of freedom), there are only $27$ control points in the global array (with $54$ degrees of freedom) due to incidence constraints, Dirichlet boundary conditions, and rigidity constraints.} 
\label{fig:numerical-methods}   
\end{figure}

\subsubsection{Isogeometric Analysis (IGA)}\label{sec:geom-rep}
IGA~\cite{HughesEtal2005} has emerged as a powerful computational method that integrates computer-aided design (CAD) and finite element analysis (FEA)~\cite{Hughes2012}, offering an exact geometric representation and seamless transition from geometric modeling to numerical simulations. 
The fundamental principle of IGA lies in utilizing smooth basis functions employed in CAD modeling, such as non-uniform rational B-splines (NURBS), to both represent the geometry and approximate the unknown fields within the domain of interest.
This avoids the need for explicit mesh generation and the associated error introduced during mesh approximation. 

In IGA the geometric domain is represented as composed of $P$ patches, each parameterized by NURBS. NURBS are built from B-splines that are defined recursively starting with piecewise constant ($p=0$) basis functions
\begin{equation}
  B_{i,0}(\xi)=\begin{cases}
    1 & \text{if $\xi_i \leq \xi < \xi_{i+1}$},\\
    0 & \text{otherwise},
  \end{cases}
\end{equation}
where $\xi_i \in \mathbb{R}$ is the $i$th knot in a knot vector $\Xi=\{\xi_1,\xi_2,\cdots,\xi_{n+p+1}\}$ that is composed of a sequence of non-decreasing real numbers in the parametric space $\xi$.
$n$ is the number of basis functions and $p$ is the polynomial order.
For $p=1,2,3,\cdots$, the basis functions are defined by
\begin{equation}
B_{i,p}(\xi)=\frac{\xi - \xi_i}{\xi_{i+p}-\xi_i}B_{i,p-1}(\xi) + \frac{\xi_{i+p+1}-\xi}{\xi_{i+p+1}-\xi_{i+1}}B_{i+1,p-1}(\xi).
\end{equation}
NURBS basis functions $N_{i,p}$ can be characterized based on B-spline basis functions by assigning a positive weight $w_i$ to each basis function 
\begin{equation}
N_i^p(\xi)=\frac{w_i B_{i,p}(\xi)}{\sum_{k=1}^n w_k B_{k,p}(\xi)}.
\end{equation}
Two dimensional (2D) NURBS basis functions with order $p$ and $q$ corresponding to knot vectors $\Xi=\{\xi_1,\xi_2,\cdots,\xi_{n+p+1}\}$ and $H=\{\eta_1,\eta_2,\cdots,\eta_{m+q+1}\}$, respectively, are constructed as tensor products
\begin{equation}
N_{i,j}^{p,q}(\xi,\eta)=N_i^p(\xi)N_j^q(\eta)=\frac{w_{i,j} B_{i,p}(\xi)C_{j,q}(\eta)}{\sum_{k=1}^n\sum_{l=1}^m w_{k,l} B_{k,p}(\xi)C_{l,q}(\eta)}.
\end{equation}
($\xi,\eta) \in [0, 1]$ corresponds to the knot span in IGA, which is often referred to as \emph{parent domain}. 
The mapping from the parent domain to the physical domain $(x,y)$ is created by a linear combination of NURBS basis functions, where the coefficients of basis functions are referred to as control points $\textbf{P}_{ij}$ (analogous to nodal coordinates in FEA). 
Therefore, the geometry of a given patch $\bm G^e(\xi)$ in two dimensions (NURBS surfaces) is represented as
\begin{equation}
\bm G^e(\xi,\eta) = \sum_{i=1}^{n}\sum_{j=1}^{m}\textbf{P}_{ij} N_{i,j}^{p,q}(\xi,\eta).
\end{equation}

In IGA, identical basis functions are used to represent the geometry and approximate the unknown field variable of a system of partial differential equations. As a result, the field variable within a patch i.e., displacement $\mathbf u^e$ can be described by
\begin{equation}
\mathbf u^e(\xi,\eta) = \sum_{i=1}^{n}\sum_{j=1}^{m}\textbf{U}_{ij} N_{i,j}^{p,q}(\xi,\eta),
\end{equation}
where $\textbf{U}_{ij}$ are control (nodal) variables, which can be computed using a Galerkin method or by finding the stationary points of a strain energy function as described in Section~\ref{sec:solver}.

We visualize the decomposition of a representative bistable unit cell into 18 patches in Figure~\ref{fig:numerical-methods}a.
All patches use the same B-spline basis functions with $n=m$ (number of basis functions or control points) and $p=q$ that are specified as parameters in our framework.
In all of our simulations, we use piecewise quadratic B-spline basis functions ($p=3$), and each control point represents two degrees of freedom in 2D space, i.e., $x_{ij}$ and $y_{ij}$.
Control points of neighboring patches have incidence constraints.
Visualization of the constraints between control points in different patches, for a representative beam composed of four patches with 16 control points in each, is shown in Figure~\ref{fig:numerical-methods}b.

\subsubsection{Material Model}\label{sec:matmodel}
We use a nearly incompressible neo-Hookean material model~\cite{Ogden1997}, which is a commonly used constitutive model that can effectively capture the mechanical behavior of elastomeric materials such as TPUs typically employed in manufacturing of flexible multistable structues.
The elastic properties of a neo-Hookean material are described by a hyperelastic strain energy density function.
This function, $W(\mathbf F)$, is independent of the path of deformation and is a function of the deformation gradient tensor, $F_{ij} = \nicefrac{\partial x_i}{\partial X_j}$ with $\mathbf X$ represents the undeformed reference configuration, and $\mathbf x$ is the deformed current configuration. In our case, $W(\mathbf F)$ is defined as

\begin{equation}
    W(\mathbf F) = \frac{\mu}{2}(I_1 - 2 - 2 \ln J) + \frac{\lambda}{2} (\ln J)^2,
\end{equation}
where $J = \det(\mathbf F)$, $I_1 = \text{tr}(\mathbf F^T \mathbf F)$, and $\mu = E/2(1+\nu)$ and $\lambda = \nu E/(1+\nu)(1-2\nu)$ are Lam\'{e} parameters of a material with Young's modulus $E$ and Poisson's ratio $\nu$. The first term is the distortional component (independent of volume change) and the second term is the near-incompressibility constraint. We can solve for the displacement $\mathbf u$ by finding the stationary point of the stored energy $\Psi(\mathbf u)$ from material mechanical deformation under external traction $\mathbf T$,
\begin{equation}\label{enrgy_eq}
    \min \bigg[\Psi(\mathbf u) = \int_{\Omega} W(\mathbf F) d\mathbf X - \int_{\partial \Omega^N} \mathbf T \cdot \mathbf u d\mathbf S\bigg],
\end{equation}
constrained by Dirichlet boundary conditions i.e., $\mathbf u = \mathbf u_b$ on $\partial \Omega^D$.
We find that the energy minimization formulation dovetails extremely well with the use of automatic differentiation.

\subsubsection{Differentiable Elasticity Solver}\label{sec:solver}
We develop an in-house end-to-end differentiable IGA framework to solve the nonlinear elasticity boundary value problem.
Our implementation is entirely in JAX~\cite{JAX}, a scientific computing and machine learning framework in Python with automatic differentiation and compilation to GPUs.

\paragraph{Brief Introduction to Automatic Differentiation}
Automatic Differentiation (AD) is a way to systematically compute mathematical derivatives of a function in a computer program. The basic idea is to write the function by composing primitive operations in an AD framework, such as JAX~\cite{JAX}. These primitives each have known pre-implemented Jacobians, and are fundamental enough (e.g.\ matrix-vector multiplication, $e^{x}$) to be composed into functions that express most operations used in scientific computing.
%written in an AD framework, such as JAX~\cite{JAX}. The basic idea is that the framework is composed of primitive operations, each of which has known pre-implemented Jacobians. These primitives can then be combined to form functions that express most operations that are used in scientific computing.
The derivatives are formed by tracing the primitive operations and combining their Jacobians according to the chain rule of calculus. When used correctly, AD computes exact derivatives and has low overhead; if a function takes time $T$ to compute, the gradient and function together can be computed in time $4T$, the \emph{Cheap Gradient Principle}~\cite{griewankwalther2008}. This constant could be less if the gradient function is precompiled, as is supported in JAX.

The fundamental operations in JAX used for differentiation are \textsc{jax.vjp} and \textsc{jax.jvp}. The former computes vector-Jacobian products with a multi-input multi-output function, and the latter computes Jacobian-vector products. Concretely, given a vector-valued function $f: \mathbb{R}^n \to \mathbb{R}^m$, the Jacobian is a matrix-valued function $\mathcal{J}[f]: \mathbb{R}^n \to \mathbb{R}^{n \times m}$. A VJP is a function of input $x \in \mathbb{R}^n$ and an \emph{adjoint} vector $v \in \mathbb{R}^m$ which computes $\left[\mathcal{J}[f](x)\right]^T v \in \mathbb{R}^n$. Note that in the case where $f$ is a scalar-output function, the VJP exactly corresponds to the gradient (when the adjoint $v = 1$). A JVP on the other hand is a function of input $x \in \mathbb{R}^n$ and a \emph{tangent} vector $v \in \mathbb{R}^n$ and computes $\left[\mathcal{J}[f](x)\right] v \in \mathbb{R}^m$. Note that in both cases, the full Jacobian is never necessarily formed; most AD frameworks will compute these as implicit matrix-vector products. They produce exact results to numerical accuracy, and are cheap to compute (comparable to runtime of original function).

We could also combine these fundamental operations to form a Hessian-vector product with the Hessian of a scalar-output function, again in time proportional to $T$, the time for function evaluation~\cite{griewankwalther2008}. We use this fact extensively in our JAX-based solver when we use second-order methods to minimize the stored energy of the material for mechanical deformation computation.

\paragraph{Numerical Forward Solve}
In our framework, we solve the boundary value problem by directly optimizing the energy as a function of displacement $\mathbf u$, subject to prescribed boundary conditions. 
We define a function in JAX to compute $\Psi(\mathbf u)$, and optimize this using a second order optimization technique, such as Newton's method.
%Since by construction $\mathbf q$ is an unconstrained parameterization, as described in Section~\ref{sec:constraints}, we are able to use standard unconstrained optimization algorithms, such as Newton's method.
The Hessian and gradient are computed using the built-in automatic differentiation engine of JAX, with special regards to the sparsity pattern of the Hessian. 

Concretely, given $\mathbf u$ representing the displacement of the body, we compute the strain energy $\Psi(\mathbf u)$ described in Equation~\ref{enrgy_eq} with primitive JAX operations. We then call the standard library function \textsc{jax.grad} (built off of \textsc{jax.vjp}) to efficiently get the gradient function with respect to $\mathbf u$. We could then apply standard gradient descent (explicit) methods to optimize $\Psi$ with only first-order information, but we find that second order (implicit) methods are needed for practical fast optimization.

Using the standard library function \textsc{jax.hessian} will compute the full dense Hessian, which would be impractical to use directly since it would require $n_{\text{dof}}$ JVPs and take $n_{\text{dof}} \times n_{\text{dof}}$ memory. Instead, we take advantage of the fact that the sparsity pattern of the Hessian is static and depends solely on the connectivity of the degrees of freedom, and use Hessian-vector products (HVPs) to probe the function to materialize the sparse Hessian.
%It is well-known in the Automatic Differentiation literature that computing Hessian-vector products costs within a constant factor of computing $\Psi(\mathbf q)$, and this is supported in JAX using the function \textsc{jax.jvp}.
In particular, we apply the Jacobian coloring algorithm~\cite{sparsehess}, which gives us $C$ special vectors $v_i \in \mathbb{R}^\text{dof}$ such that the HVPs with those vectors give us all the entries of the sparse Hessian. We then reshape those entries into a CSR sparse matrix, and use GMRES with an iLU preconditioner to solve the linear systems. For a better conditioned optimization, we use numerical continuation via load incrementing~\cite{bonet_wood_2008}.

\paragraph{Numerical Backward Solve (Adjoint Optimization)}
We leverage IGA to represent both our geometry and solution in the same B-spline basis.
This allows us to directly construct a differentiable map from geometry parameters to solution basis, and by using adjoint methods a differentiable map from geometry parameters to a specified loss function, enabling utilization of gradient-based optimization methods to find geometry parameters minimizing an arbitrary loss function.
This differs from classic shape optimization approaches such as that in~\cite{dokken2020automatic} as our approach will not require black-box meshing; the mapping from geometry to solution basis is well-defined and fully differentiable.

To differentiate with respect to the geometry parameters and material properties, we use the adjoint method. 
We wish to optimize a loss function $\mathcal{L}(\mathbf u, \theta)$ by varying $\theta$ (representing geometry and material parameters). 
In addition to possible explicit dependence of $\mathcal{L}$ on $\theta$, $\mathbf u$ is an implicit function of $\theta$ through the condition $\nabla \Psi(\mathbf u; \theta) = 0$, since our solution of the forward problem is a stationary point of $\Psi$. 
For simplicity, say $\theta \in \mathbb{R}^m$ and $\mathbf u \in \mathbb{R}^n$.

We wish to compute the total derivative $\nicefrac{d\mathcal{L}}{d\theta}$. Using the chain rule
\begin{equation}
\frac{d\mathcal{L}}{d\theta} = \frac{\partial \mathcal{L}}{\partial \mathbf u} \frac{d\mathbf u }{d\theta} + \frac{\partial \mathcal{L}}{\partial \mathbf \theta},
\end{equation}
where the terms $\nicefrac{\partial \mathcal{L}}{\partial \mathbf u} \in \mathbb{R}^n$ and $\nicefrac{\partial \mathcal{L}}{\partial \mathbf \theta} \in \mathbb{R}^m$ can be computed efficiently using the automatic differentiation capabilities in JAX. Meanwhile, $\nicefrac{d\mathbf u }{d\theta} \in \mathbb{R}^{n \times m}$ is a potentially large matrix that we do not need to materialize; we are only interested in the product $\left[\nicefrac{d\mathbf u}{d\theta}\right]^T \nicefrac{\partial \mathcal{L}}{\partial \mathbf u}$. To compute this, note that we can use implicit differentiation on the stationarity condition, $\nabla \Psi(\mathbf u; \theta) = 0$:
\begin{equation}\label{adjoint}
 \frac{\partial \nabla \Psi}{ \partial \mathbf u} \frac{d\mathbf u}{d\theta} + \frac{\partial \Psi}{ \partial \theta} = 0 \implies 
\frac{d\mathbf u}{d\theta} = - \left[\frac{\partial \nabla \Psi}{ \partial \mathbf u}\right]^{-1} \frac{\partial \Psi}{ \partial \theta}.
\end{equation}
Left multiplying Equation~\ref{adjoint} by $\nicefrac{\partial \mathcal{L}}{\partial \mathbf u}$ gives
\begin{equation}
\frac{\partial \mathcal{L}}{\partial \mathbf u} \frac{d\mathbf u}{d\theta} = - \frac{\partial \mathcal{L}}{\partial \mathbf u} \left[\frac{\partial \nabla \Psi}{ \partial \mathbf u}\right]^{-1} \frac{\partial \Psi}{ \partial \theta}.
\end{equation}
Defining $\lambda := \nicefrac{\partial \mathcal{L}}{\partial \mathbf u} \left[\nicefrac{\partial \nabla \Psi}{ \partial \mathbf u}\right]^{-1} $, we now have an algorithm to compute the product $\left[\nicefrac{d\mathbf u }{d\theta}\right]^T \nicefrac{\partial \mathcal{L}}{\partial \mathbf u}$:
\begin{enumerate}
    \item Compute $\lambda$ as the solution to the linear system
    \begin{equation}
    \left[\frac{\partial \nabla \Psi}{ \partial \mathbf u}\right]^T \lambda^T = \left[\frac{\partial \mathcal{L}}{\partial \mathbf u}\right]^T.
    \end{equation}
    Note that the matrix on the LHS is exactly the Hessian of the energy at the solution point, and since the Hessian is symmetric this linear system can be solved in the same way as we did for the forward problem.
    \item Compute the answer
    \begin{equation}
    \frac{\partial \mathcal{L}}{\partial \mathbf u} \frac{d\mathbf u}{d\theta} = - \lambda \frac{\partial \Psi}{ \partial \theta}.
    \end{equation}
\end{enumerate}
Now that we have the total derivative (gradient) of the loss function $\mathcal{L}$ with respect to the geometric and material parameters, we can use standard first order optimization to minimize the loss function. 
We impose constraints on geometric and material parameters by using projected gradient descent; after each gradient update, we project into the set defined by our constraints. 
In our application, we only have simple box constraints, which amounts to clipping the parameters after each update.

\begin{figure}[h!] 
\centering
\includegraphics[width=1.0\textwidth]{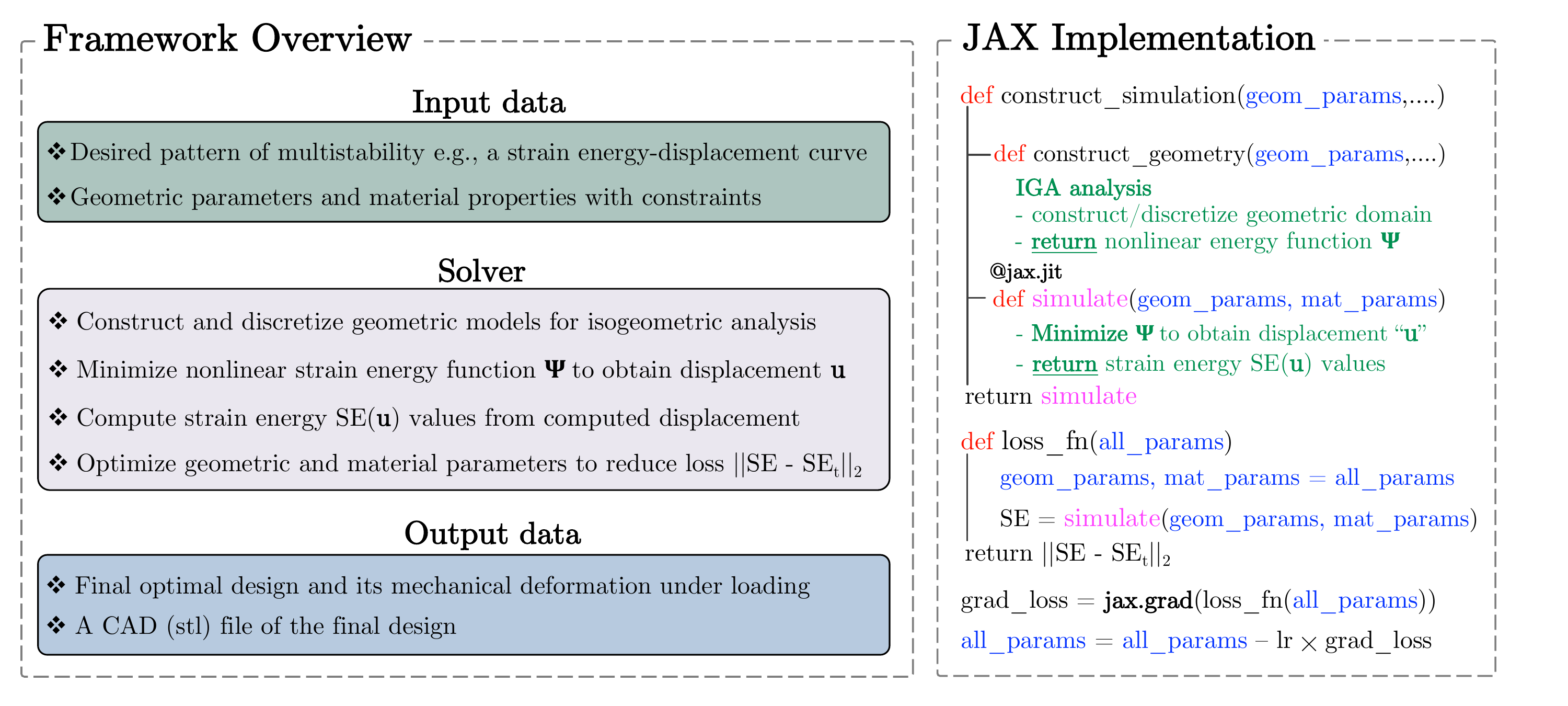}
\vspace{-7mm}
\caption{An overview of the modeling steps of our automated framework implemented in JAX for multistability design.} 
\label{fig:overview-jax}
\end{figure}

The overview and JAX implementation of our computational framework is shown in Figure~\ref{fig:overview-jax}. 
This entailed giving a desired pattern of multistability e.g., a curve of mechanical strain energy-displacement, and geometrical and material properties with constraints as the only required inputs to the framework. 
Given geometric and material design parameters, our solver simulates the mechanics of the structure by minimizing the nonlinear strain energy density function (forward solve), computes the gradient with respect to input design parameters and optimizes them via gradient-based optimization techniques (backward solve) to achieve the target multistability behavior. 
Outputs include a final optimal design, its mechanical responses under uniaxial compression loading, and an stl file of the model for 3D printing. 

\section{Results}\label{sec:Results}
Here we apply the proposed computational framework to design several bistable and multistable structures and compare the results to those derived from experiment. To ensure consistent comparisons, the same loading and boundary conditions, and material properties are used in simulations to appropriately resemble the experiments. Multiple tests were conducted for each sample to get average experimental load-displacement curves, the integral of which give strain energy values using  Simpson's rule. In all simulations, the geometric parameters and the Young's modulus are optimized to reach a loss value less than 0.01. Mesh independence analysis was performed to ensure simulation results were not affected by insufficient mesh resolution, where generally seven control points for each patch was observed to be enough. In this study, all numerical results are obtained by applying a uniform compression displacement loading at the top edge of the structure. 
%\todo{RPA: this last sentence is a mess.} We note that although our framework can design multi-column SIMS, here we mainly focus on \textit{single}-column structures. The main reason is that multi-column SIMS have the same stability locations as the corresponding single-column structure with identical material and mechanical properties, and their strain energies are equal to values from a single-column model times the number of columns.\todo{RPA: assuming identical columns an uniaxial loads...} Therefore, one reasonable strategy could be dealing with normalized strain energies, however, actual values are presented in this work. \todo{RPA: reasonable strategy for what?}

\subsection{Bistablity analysis of a unit cell}\label{sub:bistablity}
\begin{figure}[h!] 
\centering
\includegraphics[width=1.0\textwidth]{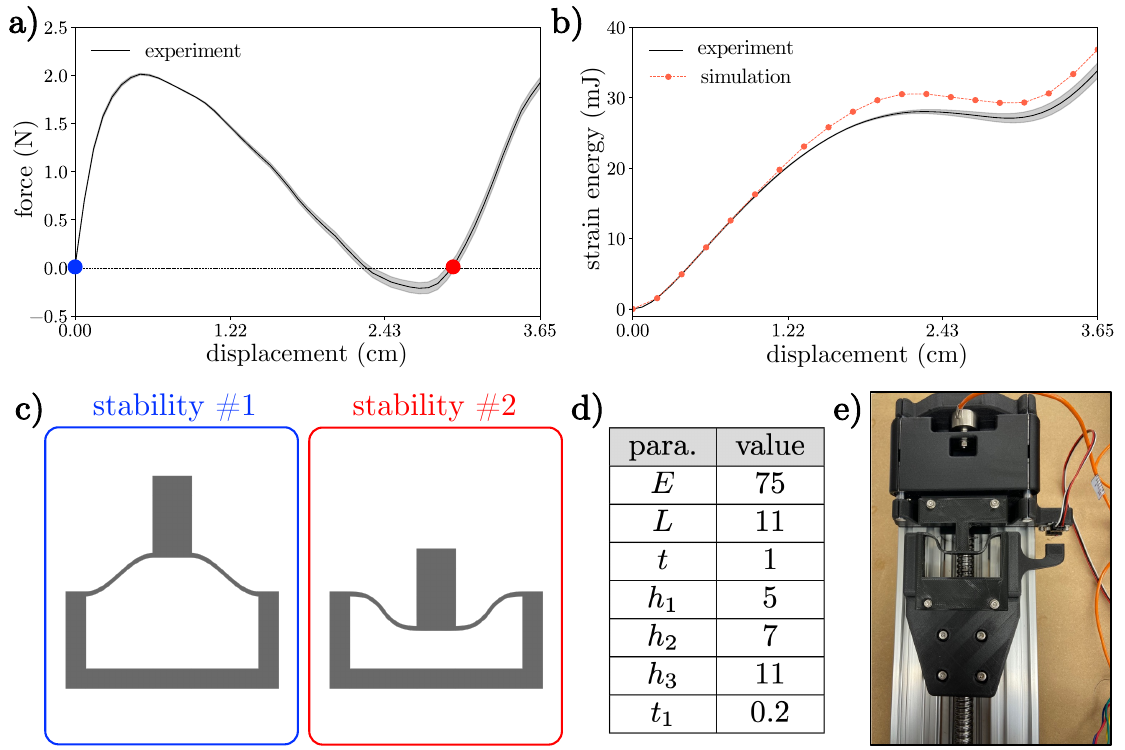}
\caption{\textbf{a)}\ Experimentally measured forces against vertical displacement for a unit cell under uniaxial loading.\ \textbf{b)}\ Elastic strain energies versus displacement from simulation and experiment.\ \textbf{c)}\ Stable configurations of the unit cell.\ \textbf{d)}\ Geometric parameters (in cm) and material property (in MPa) of the model.\ \textbf{e)}\ Our experimental setup to perform bistability test.} 
\label{fig:ms_f2}   
\end{figure}

Figure~\ref{fig:ms_f2}a shows experimentally measured forces against vertical displacement for a unit cell with geometric and material properties given in Figure~\ref{fig:ms_f2}d. The bistability deformation behavior is predicted by the experiment. The force initially increases along with the vertical displacement until a maximum critical compression force is reached, from which point snap-through instability of the structure is triggered and negative stiffness is observed. The force reduces to zero, and then becomes negative representing tension. 

The elastic strain energy versus displacement curves are generally utilized for stability analysis of multistable structures, where the global and local minima of the strain energy curve are indicative of stable states of a structure. The variation of strain energy during the loading process is demonstrated in Figure~\ref{fig:ms_f2}b from experiments and simulation. The strain energy values are the areas under the corresponding force-displacement curve from the experiment, and are computed in our isogeometric analysis. There is a global minimum at $x=0$ and a local minimum at $x=2.71$ in Figure~\ref{fig:ms_f2}b, representing two stable positions of the single unit cell shown in Figure~\ref{fig:ms_f2}c that are identically predicted by both the experiment and simulation. The local maximum energy is an unstable equilibrium position that defines the amount of energy barrier to overcome to switch from first stability to the second. Once the hill of the barrier is passed, the structure automatically jumps into a lower-energy stable configuration by releasing the additional energy.

The quantitative discrepancies between experiment and simulation results generally stem from approximated constitutive modeling assumptions, material properties estimation, and difficulties in experimentally reproducing exact boundary conditions used in simulations.

\subsection{Geometric optimization for bistability design}\label{sub:opt-su}
Here, we study the ability of our framework to optimize geometric parameters (while assuming a constant material properties i.e., $E=70$MPa) to design a single unit cell with a desired pattern of bistability determined by a strain energy-displacement curve, which is obtained from simulating uniaxial compression loading of a unit cell with geometric parameters given in Table~\ref{table:ms_t1}. 

\begin{figure}[h!] 
\centering
\includegraphics[width=1.0\textwidth]{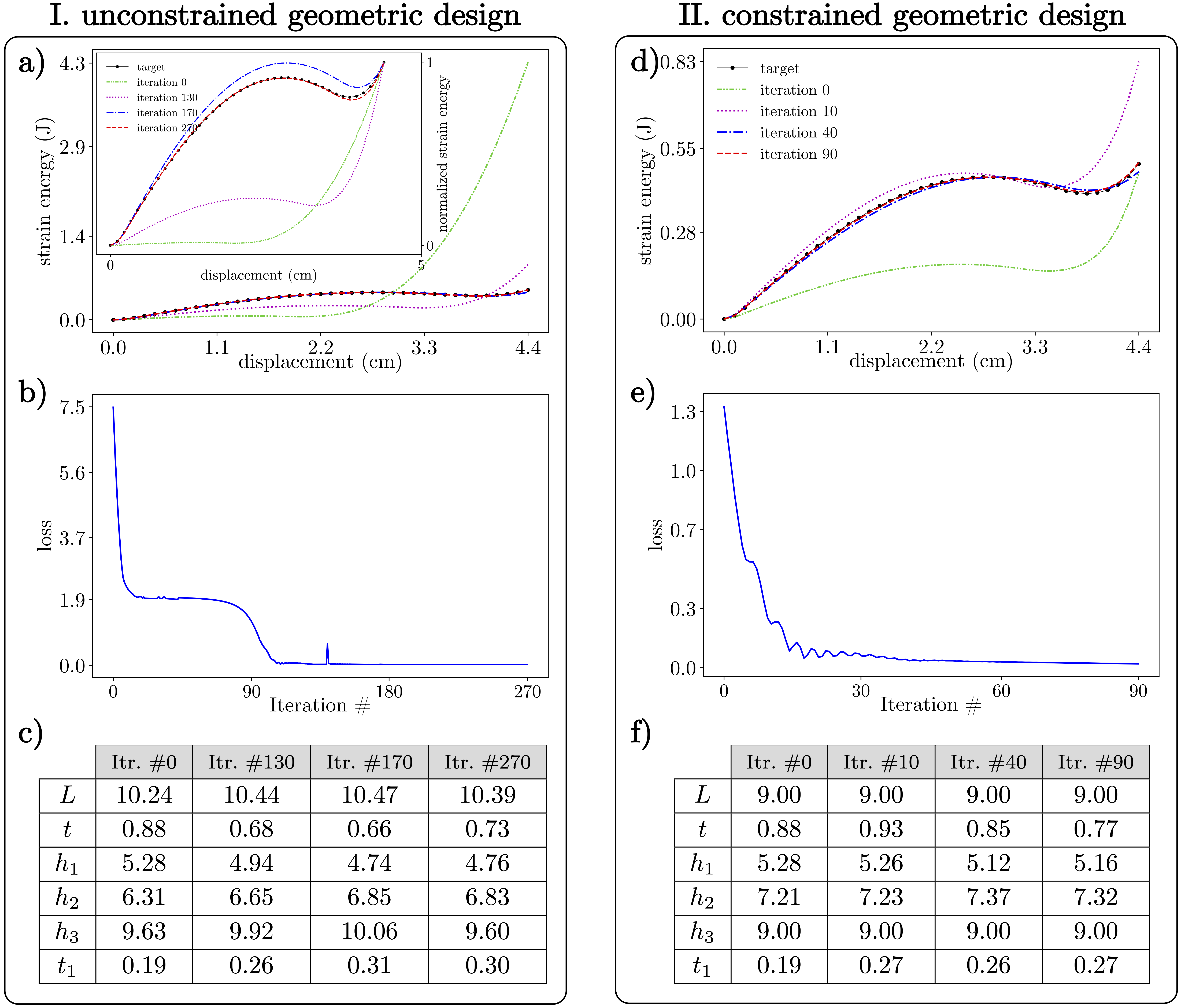}
\caption{\textbf{a,d)}\ Strain energy-displacement curves against displacement from uniaxial loading simulations of a designed unit cell with geometric parameters at different optimization iterations.\ \textbf{b,e)}\ Loss values versus iteration numbers.\ \textbf{c,f)}\ Geometric properties of the model (in cm) at different iterations.} 
\label{fig:ms_f3}   
\end{figure}

To evaluate the robustness of our framework, we tested two different strategies: \textbf{I)} Initial geometric parameters of the design structure (cf.\ Iteration $\#0$ at Figure~\ref{fig:ms_f3}c) are selected such that the initial strain energy-displacement curve remarkably deviates from the target curve as shown in Figure~\ref{fig:ms_f3}a. Even with such a substantial discrepancy between initial and target strain energy values, our framework was able to find a bistable structure with a strain energy-displacement curve in agreement with the target curve after a reasonable number of iterations.  \textbf{II)} Geometric design constraints are used while matching a target strain energy-displacement curve. For example, here we used same geometric parameters as in strategy \textbf{(I)} except $h_2$ is carefully chosen to make the stability location of the initial design close to the desired value while we are interested in a unit cell with $L=9$ and $h_3=9$ (Figure~\ref{fig:ms_f3}f). Our framework was able to design a bistable structure with these geometric constraints that has desired mechanical responses under the applied compression loading condition.

\begin{table}[h!]
\begin{center}
\begin{tabular}{|c|c|c|c|c|c|c|}
\rowcolor{gainsboro}
\hline
\cellcolor{white} parameter  &  $L$ &  $t$ &  $h_1$ &  $h_2$ &  $h_3$ &  $t_1$\\ 
\hline
value & 11.34 & 1.24 & 4.15 & 6.28 & 10.17 & 0.28\\ 
\hline
\end{tabular}
\end{center}
\vspace{-4mm}
\caption{Geometric parameters of a single unit cell model (in cm) used to create the target strain energy-displacement curve in Figure~\ref{fig:ms_f3}.}
\label{table:ms_t1}
\end{table}

\subsection{Multistability analysis and SIMS design}\label{sub:multistablity}
\begin{figure}[h!] 
\centering
\includegraphics[width=1.0\textwidth]{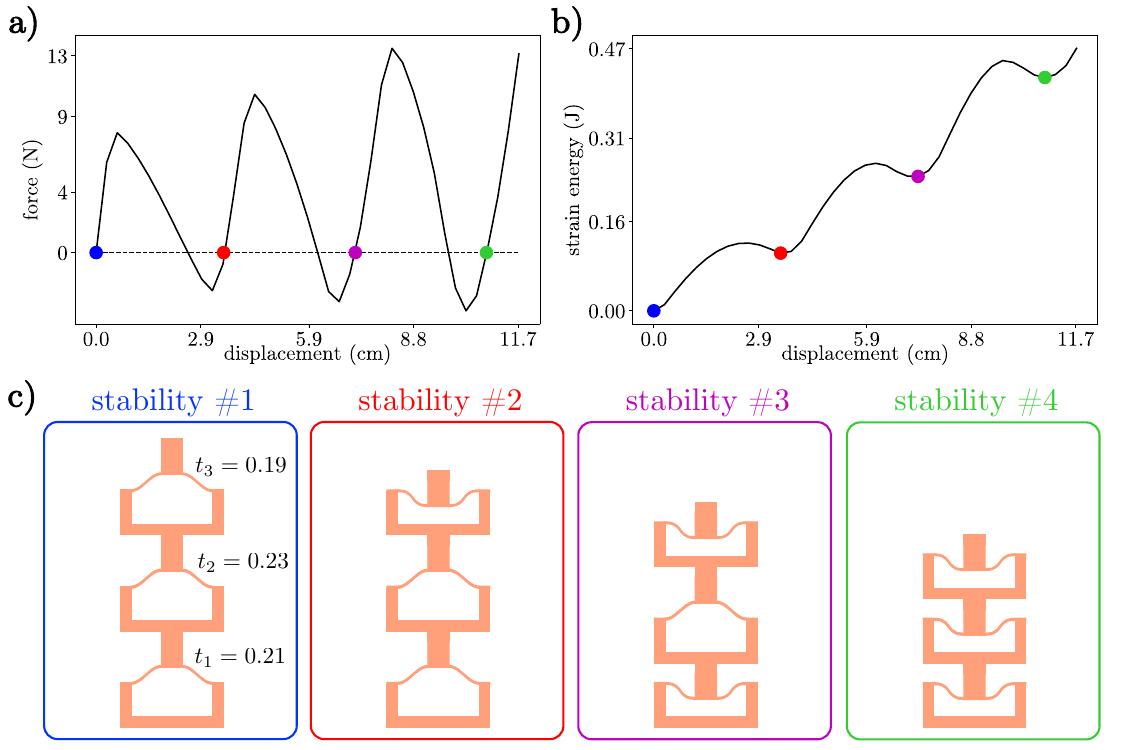}
\caption{\textbf{a)}\ Computed forces against vertical displacement from uniaxial loading of a single column SIMS composed of three identical unit cells with different in-plane beam thicknesses. \textbf{b)}\ Elastic strain energies versus displacement from simulation. \textbf{c)}\ Four stable configurations of the single column SIMS.} 
\label{fig:ms_f4}   
\end{figure}

Multistability is studied here for a SIMS that includes three identical unit cells in series with varying in-plane beam thicknesses. Figure~\ref{fig:ms_f4}a illustrates negative stiffness and tensile forces observed three times during loading process in the simulation. The presence of a global minimum, representing the initial undeformed configuration, and three local minima in the computed strain energy landscape shown in Figure~\ref{fig:ms_f4}b indicates the multistability of this structure. The sequence of deformation the structure switches between stable configurations under uniaxial loading predicted by simulation is consistent with previous studies proposing the rows collapse in the order that thickness values increase~\cite{CheEtal2016}. Therefore, according to thickness values in Table~\ref{table:ms_t2} (the unit cells are ordered from bottom to top), the second stability corresponds to when only the top cell was collapsed followed by the third stable configuration that is achieved when top and bottom cells were collapsed (Figure~\ref{fig:ms_f4}c). Finally, the last stability happened by applying $\sim13.5$ (N) compression force required for collapsing all three cells. 

\begin{table}[h!]
\begin{center}
\begin{tabular}{|c|c|c|c|c|c|c|c|c|}
\rowcolor{gainsboro}
\hline
\cellcolor{white} parameter  &  $L$ &  $t$ &  $h_1$ &  $h_2$ &  $h_3$ &  $t_1$ &  $t_2$ &  $t_3$\\ 
\hline
 value & 12.21 & 1.25 & 5.32 & 7.24 & 11.45 & 0.21 & 0.23 & 0.19\\ 
\hline
\end{tabular}
\end{center}
\vspace{-4mm}
\caption{Geometric parameters of a single column SIMS model (in cm) used to perform computational multistability analysis shown in Figure~\ref{fig:ms_f4}, and to create the target strain energy-displacement curve in Figure~\ref{fig:ms_f5}.}
\label{table:ms_t2}
\end{table}

\begin{figure}[h!] 
\centering
\includegraphics[width=1.0\textwidth]{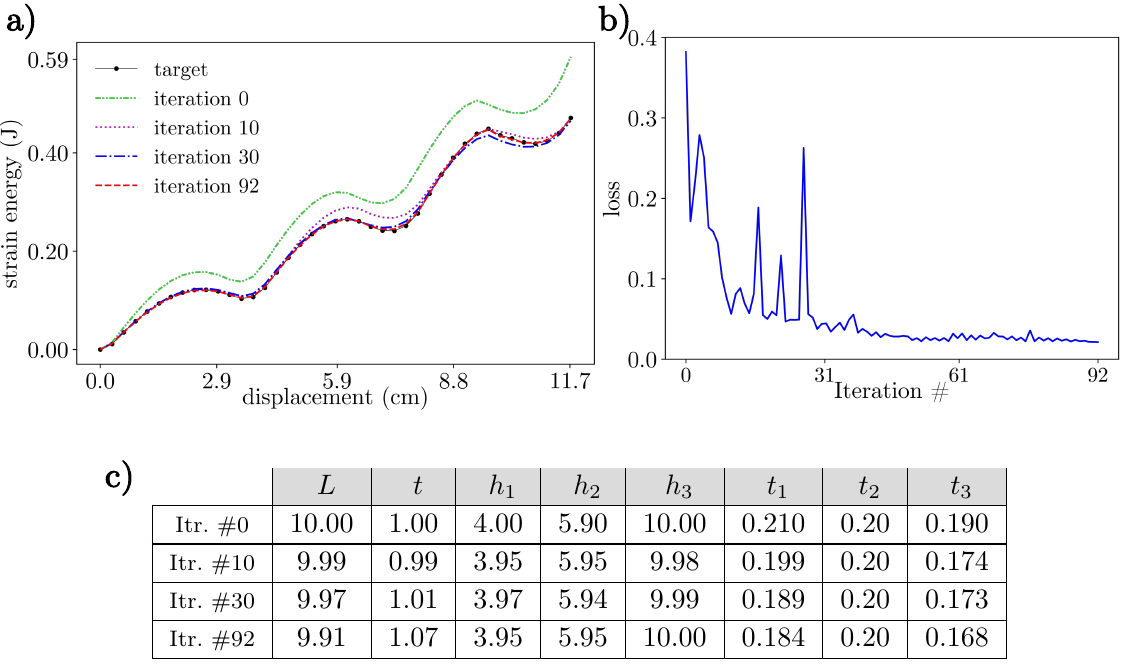}
\caption{\textbf{a)}\ Strain energy-displacement curves against displacement from uniaxial loading simulations of SIMS with geometric parameters at different optimization iterations.\ \textbf{b)}\ Loss values versus iteration numbers.\ \textbf{c)}\ Geometric properties of the model (in cm) at different iterations.} 
\label{fig:ms_f5}   
\end{figure}

The strain energy curve in Figure~\ref{fig:ms_f4}b is now used as a desired pattern of multistability in designing a new SIMS with $L\leq10$, $h_3\leq10$, and $t_i\leq0.2$ (cf.\ Table~\ref{table:ms_t2}). The eight initial geometric parameters of the model are given in Figure~\ref{fig:ms_f5}c. Our framework optimized these geometric parameters in 92 iterations to design a structure that satisfies given geometric design constraints, and has the target strain energy-displacement curve (Figure~\ref{fig:ms_f5}a).

\begin{figure}[h!] 
\centering
\includegraphics[width=1.0\textwidth]{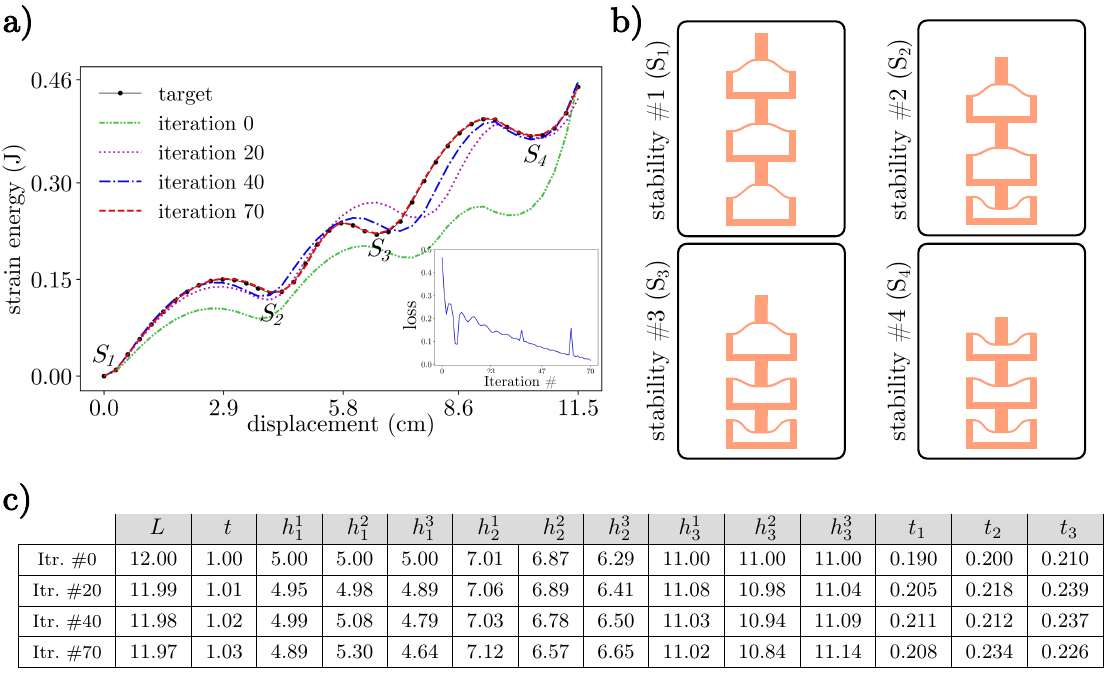}
\caption{\textbf{a)}\ Strain energy-displacement curves against displacement from uniaxial loading simulations of SIMS with geometric parameters at different optimization iterations, and loss values versus iteration numbers.\ \textbf{b)}\ Four stable configurations of the final designed SIMS.\ \textbf{c)}\ Geometric properties of the model (in cm) at different iterations.} 
\label{fig:ms_f6}   
\end{figure}

We also evaluated the utility of our computational framework for designing a more complex SIMS with three non-identical unit cells resulting in 14 geometric design parameters to optimize. The desired multistablity pattern in Figure~\ref{fig:ms_f6}a includes four stable locations at $x=0, 4.03, 6.61, 10.35$ (cm), and energy barriers to switch from $S_1 \to S_2$, $S_2 \to S_3$, and $S_3 \to S_4$ are 0.15, 0.11, and 0.18 (J), respectively. Figure~\ref{fig:ms_f6}b demonstrates the deformation sequence of multistability pattern for the final design with given geometric parameters in Figure~\ref{fig:ms_f6}c, i.e., the bottom row collapses first, then the middle row and finally the top row.        

\begin{figure}[h!] 
\centering
\includegraphics[width=1.0\textwidth]{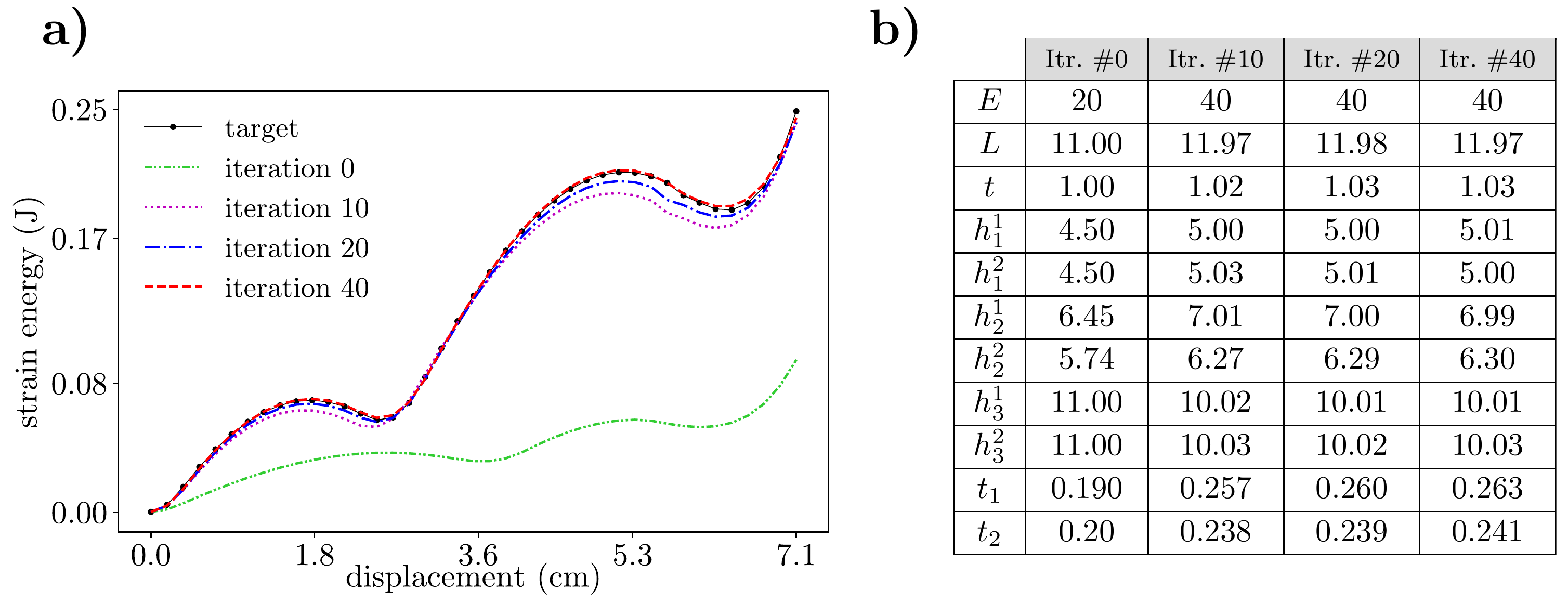}
\caption{\textbf{a)}\ Strain energy-displacement curves against displacement from uniaxial loading simulations of a SIMS with geometric parameters at different optimization iterations.\ \textbf{b)}\ Geometric properties (in cm) and Young's modulus (in MPa) of the model at different iterations.} 
\label{fig:ms_f7}   
\end{figure}

Here, our framework is utilized to simultaneously optimize geometric parameters and Young's modulus of a material to build a SIMS with particular multistability behavior. This enables designing a desired multistability pattern with specific material choices. We examined this capability of our framework for a SIMS with two non-identical unit cells to achieve a target multistability design in Figure~\ref{fig:ms_f7}a. The geometric parameters and material property used to create the target strain energy-displacement curve are presented in the first row of Table~\ref{table:ms_t3}. We assumed that we have three material choices with $E=20, 40$, and $80$ MPa and no constraints on geometry in our design. The framework proposed a structure that has the desired multistability behavior with geometric parameters given in Figure~\ref{fig:ms_f7}b using a material with $E=40$ MPa, where the target strain energy-displacement curve was created assuming an $E=60$ MPa. 

In another experiment, we investigated the possibility of creating the same target strain energy curve in Figure~\ref{fig:ms_f7}a using materials with Young's modulus of 20 MPa and 100 MPa and initial geometric parameters in the first column of Figure~\ref{fig:ms_f7}b. The geometric parameters of the final design with identical strain energy curves for each material are presented in Table~\ref{table:ms_t3}.          

\begin{table}[h!]
\vspace{4mm}
\begin{center}
\begin{tabular}{ccccccccccc}
\hline
\multicolumn{1}{c}{$^*$target curve} & \multicolumn{10}{c}{\cellcolor{gainsboro}{geometric parameters}} \\
\cmidrule(rl){2-11}
\rowcolor{gainsboro}
{Young's Modulus ($E$)} & {$L$} & {$t$} & {$h_1^1$} & {$h_1^2$} & {$h_2^1$} & {$h_2^2$} & {$h_3^1$} & {$h_3^2$} & {$t_1$} & {$t_2$}\\
\hline
\cellcolor{EB} \hspace{2mm}60$^*$ &  12.00 &  1.00 &  5.00 &  5.00 &  7.00 &  6.30 &  10.00 &  10.00 &  0.230 &  0.210 \\ 
\cellcolor{EB} 40 &  11.97 &  1.03 &  5.01 &  5.00 &  6.99 &  6.30 &  10.01 &  10.03 &  0.263 &  0.241 \\ 
\cellcolor{EB} 20 &  11.88 &  1.11 &  5.02 &  5.01 &  6.99 &  6.29 &  9.89 &  10.44 &  0.324 &  0.298\\ 
\cellcolor{EB} 100 &  12.05 &  0.95 &  4.99 &  5.00 &  7.01 &  6.30 &  9.94 &  9.94 &  0.196 &  0.181 \\ 
\hline
\end{tabular}
\end{center}
\vspace{-4mm}
\caption{Young's modulus (in MPa) and geometric parameters (in cm) of four different SIMS under uniaxial loading with identical strain energy-displacement curve responses shown as the target curve in Figure~\ref{fig:ms_f7}a.}
\label{table:ms_t3}
\end{table}

\section{Discussion}\label{sec:Discussion}
We have presented a framework implemented in JAX to rapidly design multistable mechanical structures. 
The proposed end-to-end differentiable IGA-based continuum mechanics solver enables geometric and material optimization of a structure to achieve a desired multistability pattern. 
Our framework can significantly facilitate and automate soft robotic design for high precision motion control (via controlling stability locations and their sequences) and optimal force/energy output (via controlling energy barriers and releases).
The inputs to the framework are geometrical and material properties, specification of boundary conditions, and a desired multistability pattern. 
The associated computational cost on a single GPU core for compilation and execution times at each optimization iteration are shown in Table~\ref{table:ms_t4} for multistable structures with different number of rows, which are obtained by simulating uniaxial compression displacement loading for 100 structural models in each case. 

\begin{table}[h!]
\begin{center}
\begin{tabular}{cccccc}
\hline
\multicolumn{1}{c}{} & \multicolumn{5}{c}{\cellcolor{gainsboro}{number of rows}} \\
\cmidrule(rl){2-6}
\rowcolor{gainsboro}
{computational time (s)} & {one} & {two} & {three} & {four} & {five} \\
\midrule
compilation & 142.90 $\pm$ 1.18 & 162.33 $\pm$ 2.32 & 191.85 $\pm$ 2.45 & 223.68 $\pm$ 3.53 & 298.56 $\pm$ 3.21 \\
execution  & 7.24 $\pm$ 1.50 & 16.39 $\pm$ 1.18 & 34.64 $\pm$ 3.26 & 51.46 $\pm$ 3.76 & 109.07 $\pm$ 5.55 \\
\hline
\end{tabular}
\end{center}
\vspace{-4mm}
\caption{Compilation and execution times (in seconds) from simulating single column SIMS with different number of rows under uniaxial compression loading (100 simulations for each).}
\label{table:ms_t4}
\end{table}

In addition to GPU acceleration that makes simulations approximately five times faster, we considered several improvement in our computational framework to speed up the simulation process. Namely, in our isogeometric analysis some of the patches within each unit cell are regarded as rigid (nine patches with star in Figure~\ref{fig:numerical-methods}a) resulting in significant reduction in number of degrees of freedom (see Table~\ref{table:ms_t5}), while strain energy values trivially change and stable locations remain the same to simulations with all deformable patches.
%($e_{SE} = \sum_{i=1}^{N}|SE_i^{R}-SE_i^{D}|/\max_n (SE_n^{D})/N$, for $N$ number of loading increment; the error is normalized by maximum values from modeling with all deformable patches, to avoid division by small values of strain energies close to zero). 
Moreover, simulations with rigid patches generally achieve the desired multistablity pattern in fewer number of optimization iterations and converge with larger displacement increment size (applied on top edge), leading to reduced total computation time as shown in Table~\ref{table:ms_t5}. 
Finally, using backtracking line search to determine appropriate step-size for Newton's method and applying an adaptive loading increment size technique enhanced both computational efficiency and automation of our framework.

\begin{table}[h!]
\begin{center}
\begin{tabular}{ccccccccc}
\hline
\multicolumn{1}{c}{} & \multicolumn{3}{c}{\cellcolor{gainsboro}{execution time (s)}} & \multicolumn{2}{c}{\cellcolor{gainsboro}{total time (s)}} & \multicolumn{2}{c}{\cellcolor{gainsboro}{number of DoF}} & \multicolumn{1}{c}{} \\
\cmidrule(rl){2-4} \cmidrule(rl){5-6} \cmidrule(rl){7-8} 
\rowcolor{gainsboro}
\# of rows & {w/o rig.} & {w/ rig.} & {w/ rig.\ (CPU)} & {w/o rig.} & {w/ rig.} & {w/o rig.} & {w/ rig.}  & error in SEs \\
\hline
1 & 14.02 & 7.77 & 43.00 & 429 & 210 & 1400 & 800 & 0.91\%  \\
2 & 41.47 & 22.02 & 104.79 & 970 & 716 & 2860 & 1601 & 1.19\% \\
3 & 92.12 & 38.99 & 195.02 & 4199 & 1475 & 4320 & 2402 & 1.37\% \\
4 & 172.27 & 65.93 & 289.17 & 29278 & 5572 & 5780 & 3203 & 1.43\%  \\
\hline
\end{tabular}
\end{center}
\vspace{-4mm}
\caption{Total simulation time and execution times at each optimization iteration (in seconds), number of degrees of freedom (DoF) for modeling with and without rigid patches, and relative error values in computed strain energies (SEs) from simulating single column SIMS with different number of rows under uniaxial compression loading by enforcing geometric constrains $L=10, h_3=10$, and $t_i \leq 0.2$. Error values between simulations with rigid patches and all deformable patches are computed as $e_{SE} = (1/N)\sum_{i=1}^{N}|SE_i^{R}-SE_i^{D}|/\max_n (SE_n^{D})$, for $N$ number of loading increment; the error is normalized by maximum values from simulations with all deformable patches to avoid division by small values of strain energies close to zero.}
\label{table:ms_t5}
\end{table}

Our framework enables reasonable flexibility in satisfying different design criteria through easily defining custom loss functions and imposing appropriate mechanical, geometrical, and material constraints. For instance, the type of multistable structures studied in this work are broadly used as soft actuators for high energy absorption and precise motion control that are merely characterized by the amount of energy barriers and stable locations, respectively, and not necessarily by the entire strain energy curve. Therefore, we can propose a custom loss function that allows to rapidly optimize the geometry of a mechanical structure to induce local extrema of a displacement strain energy curve with respect to a reference point, which determines the amount of energy barriers, energy releases, and axial displacement achieved by switching from one stable configuration to another. Namely, instead of a strain energy curve, the inputs to the framework could be desired extrema values of strain energy and the location of its minima. This crucially gives more freedom in designing multistable strucures for such applications. For example, Table~\ref{table:ms_t6} indicates strain energy values and the stablity location of a multistable structure that satisfies the designed criteria to absorb 0.12 and 0.14 (J) energy when consecutively moving from stable locations at $x=0.0, 3.83, 7.65$ (cm).     

\begin{table}[h!]
\begin{center}
\begin{tabular}{ccccccccc}
\hline
\rowcolor{gainsboro}
\multicolumn{2}{c}{maximum SEs (J)} & \multicolumn{2}{c}{minimum SEs (J)} & \multicolumn{2}{c}{energy barriers (J)} & \multicolumn{3}{c}{stability loc. (cm)} \\
\cmidrule(rl){1-2} \cmidrule(rl){3-4} \cmidrule(rl){5-6} \cmidrule(rl){7-9}
{0.118} & {0.240} & {0.101} & {0.219} & {0.118} & {0.139} & {0.0} & {3.83} & {7.65} \\
\hline
\end{tabular}
\end{center}
\vspace{-4mm}
\caption{Energy extrema and stability locations of a designed multistable structure with target energy barriers of 0.12 and 0.14 (J), and stablility at $x = 0.0, 3.83, 7.65$ (cm).}
\label{table:ms_t6}
\end{table}

Although here we focused on a particular snap-through induced multistable structure, our framework can readily be adapted for other unit cell selections with appropriate geometric parameterization such as those include stepped beams~\cite{HaghpanahEtal2016}, two camber beams~\cite{CorreaEtal2015}, twisting structures~\cite{ArrietaEtal2018} among others. However, one of the limitations is that the ability of the framework to create a multistablity pattern is restricted to the space of multistability that is achievable by that particular unit cell geometry. Topology optimization techniques may potentially mitigate such challenges and enable a broader range of multistability patterns design. A limitation in terms of numerical methods is that our current solver would have trouble at large numbers of degrees of freedom, especially when modeling 3D structures, due to the iLU factorization we perform to precondition the iterative solver. Future extensions include modern multigrid methods to scale our solver up when we consider 3D problems or 2D problems with very large number of degrees of freedom (100,000-1,000,000).

\section{Acknowledgements}
We would like to thank E Medina and PT Brun for assisting with experiments. This work was partially supported by NSF grants IIS-2007278 and OAC-2118201, the NSF under grant number 2127309 to the Computing Research Association for the CIFellows 2021 Project, and a Siemens PhD fellowship.

\section{Conflict of interest}
The authors declare no conflict of interest.

%section{Data availability}
%Our implementation is available at \href{https://github.com/PrincetonLIPS/Varmint/blob/main/projects/multistability}{\textcolor{blue}{multistability source code}}.

\clearpage
\bibliographystyle{abbrv}
\bibliography{multistability.bib}
\clearpage
\end{document}